\documentclass[12pt]{article}

\usepackage{a4}
\usepackage{parskip}

\usepackage{amsmath}
\usepackage{amsfonts}
\usepackage{amssymb}
\usepackage{theorem}
\usepackage{latexsym}

\usepackage{tikz-cd} % \usepackage{pb-diagram,lamsarrow,pb-lams}

\usepackage{epsfig}
\usepackage{graphics}
\usepackage{color}
\newtheorem{thm}{Theorem}[section]
\newtheorem{prop}[thm]{Proposition}
\newtheorem{cor}[thm]{Corollary}
\newtheorem{lem}[thm]{Lemma}
\newtheorem{expl}[thm]{Example}
\newtheorem{defi}[thm]{Definition}
\newtheorem{rem}[thm]{Remark}

\def\D{{\Bbb D}}

\def\N{{\Bbb N}}
\def\P{{\Bbb P}}

\def\R{{\Bbb R}}
\def\C{{\Bbb C}}

\def\T{{\Bbb T}}

\def\cC{{\cal C}}
\def\cD{{\cal D}}

\def\cO{{\cal O}}

\def\dim{\hbox{dim}}

\def\n*{\N_*^{m}}
\def\1{{\text{\bf 1}}}

\def\n{\vert\vert}

\def\qed{$\Box$}

\begin{document}

\title{\bf On the Hartogs extension theorem
for unbounded domains in $\mathbb C^n$}

\author{\bf Al Boggess, Roman Dwilewicz, Egmont Porten}

\date{ } %\today}

%\marginparsep=18pt
%\marginparpush=12pt
\marginparwidth=1truein
\def\margin#1{%
\vadjust{\vbox to -2.5pt
{\vss \hbox to 10pt {\hss #1\hskip.25in}}}}

\maketitle

\begin{quote} \emph{This article is dedicated to the memory of our dear colleague and friend, 
Roman Dwilewicz, who passed away on July 29, 2016. He touched many lives through 
his mathematical collaborations, his conference organizational skills, his extensive travels, 
and his ability to befriend just about anyone he met. We miss him dearly.} \end{quote}

\quad

\begin{abstract}
\noindent
Let $\Omega\subset\C^n$, $n\geq 2$, be a domain with smooth connected boundary.
If $\Omega$ is relatively compact, the Hartogs-Bochner theorem ensures that
every CR distribution on $\partial\Omega$ has a holomorphic extension to
$\Omega$. For unbounded domains this extension property may fail, for
example if $\Omega$ contains a complex hypersurface.
The main result in this paper
tells that the extension property holds if and only if the envelope
of holomorphy of $\C^n\backslash\overline{\Omega}$ is $\C^n$.
It seems that it is a first result in the literature
which gives a geometric characterization of unbounded domains in $\mathbb C^n$
for which the Hartogs phenomenon holds.
Comparing this to earlier work by the first two authors and Z.~S{\l}odkowski,
one observes that the extension problem sensitively depends on
a finer geometry of the contact of a complex hypersurface and the
boundary of the domain.
\end{abstract}

\section{Introduction}
Throughout this article we consider a domain
$\Omega=\Omega^-\subset\C^n$, $n\geq 2$,
with $\cC^\infty$-smooth connected boundary $M$.
If $\Omega$ is relatively compact in $\C^n$, the classical
Hartogs-Bochner theorem tells that every CR function on $M$
admits holomorphic extension to $\Omega$.
Via a convenient notion of weak boundary values,
this result naturally generalizes to CR distributions.

The classical Hartogs extension theorem made an important influence not
only on Complex Analysis, but also on other areas of mathematics,
like Algebraic Geometry or Partial Differential Equations.
The theorem still inspires researchers and there is a renewed interest
in recent years: Harz-Shcherbina-Tomassini \cite{HST1, HST2},
\O vrelid-Vassiliadou \cite{OV}, Damiano-Struppa-A.Vajiac-M.Vajiac
\cite{DSVV}, Pala\-mo\-dov \cite{Pal}, Ohsawa \cite{O}, Coltoiu-Ruppenthal \cite{CR},
Lewandowski \cite{lew},
and papers by the authors with other colleagues
\cite{BDS1,BDS2,BDS3,BDS4}, \cite{MP, MP2}.

A good deal of the mentioned contributions consider extension from
boundaries of unbounded domains.
Easy examples show that the Hartogs-Bochner theorem
may fail for unbounded domains, leading to the problem
to understand the precise nature of the obstacles.
The essence of the present article is a geometric characterization
of the Hartogs extension property for CR-distributions.

%\textcolor{red}{
%I would cut the following sentence, since I believe that it comes to soon
%to be understandable. In addition the matter is taken up around
%Corollary 1.5:
%Relating this to work of the first two authors with
%Z.~S{\l}odkowski, we will see that the answer to the general extension problem
%sensitively depends on a finer geometry of the contact of a complex
%hypersurface and the boundary of the domain.}

%The connectedness of $M$ yields that $\Omega^+=\C^n\backslash\overline{\Omega}$
%also is connected.

Let $S$ be a smooth real hypersurface of $\C^n$ and $\omega\subset\C^n$
a domain such that $\omega\backslash S$ has two connected components
$\omega^-$ and $\omega^+$.
A function $f\in\cO(\omega^-)$ is said to have {\it polynomial growth
at $p\in S\cap\omega$}, if there are $k\geq 0$ and
$\epsilon>0$
%and an open $\epsilon$-neighborhood $U$ of $p$
such that
\begin{equation}\label{polgrowth}
|f(z)|\,\leq \, C \,\mbox{\rm dist}(z,S)^{-k}
\end{equation}
holds for $z\in\omega^-\cap B_\epsilon(p)$.
%(where $B(p;\epsilon)$ denotes that euclidean ball centered at $p$ with radius $\epsilon$)
We say that $f$ has {\it polynomial growth towards $S$} if it has polynomial growth
at every $p\in S\cap\omega$.
It is well-known that such functions

have unique weak boundary values in $\cD'_{CR}(S)$,
the space of CR distributions on $S$,
see \cite[Ch.~VII]{BER} and also Section 2.

\begin{defi}
We say that \emph{Hartogs extension} holds for $\Omega$,
if every $u\in\cD'_{CR}(M)$ is the boundary value of
some $f\in\cO(\Omega)$ with polynomial growth along $M$.
\end{defi}

%Our goal is to study whether the Hartogs extension holds for unbounded
%domains $\Omega$.
%\begin{figure}[h]
%\parbox{.5\textwidth}{\resizebox{2.4in}{2.4in}{\input{1-omega-omega+.pstex_t}}}
%\parbox{.48\textwidth}
%{\caption{\sl Domain $\Omega$ and its complement $\Omega^+$}
%\label{fig:1-omega-omega+}}
%\end{figure}
%By the classical Hartogs-Bochner theorem, Hartogs extension always holds
%for $\Omega$ bounded.
%Moreover, some of the known proofs reduce the
%argument to the geometric fact that the envelope of holomorphy of
%$\C^n\setminus\overline{\Omega}$ is $\C^n$.
The most straightforward examples for domains without Hartogs
extension are domains containing a complex hypersurface,
but these are very far from exhausting all possible obstructions.
Despite of considerable recent activity, see
\cite{MP2, O, OV, P1, DSVV},  or older \cite{L}, to mention a few,
a satisfying understanding of Hartogs extension for unbounded domains
seems still to be missing,
even in the case that $M$ is strictly pseudoconvex at every point.
The main result of the present note is a geometric characterization that
establishes a close link to envelopes of holomorphy.

\begin{thm}\label{main} {\bf (Main Theorem)}
Let $\Omega\subset\C^n$, $n\geq 2$, be a domain with connected smooth
boundary $M$. Then Hartogs extension holds for $\Omega$ if and only
if the envelope of holomorphy of the outer domain
$\Omega^+=\C^n\setminus\overline{\Omega}$ is $\C^n$.
\end{thm}

Actually, Theorem \ref{main} is the global version
of the more general Theorem \ref{local},
where $\Omega^+$ is replaced by arbitrary outer collars attached to $M$.
Moreover, Theorem \ref{main} straightforwardly generalizes to domains in Stein manifolds.

\begin{thm}\label{main-stein} {\bf (cf.~Theorem \ref{stein})}
Let $X$ be a Stein manifold and let $\Omega\subset X$, ${\rm dim}_{\mathbb C} X\geq 2$,
be a domain with connected smooth boundary $M$.
Then Hartogs extension holds for $\Omega$ if and only
if the envelope of holomorphy of the outer domain
$\Omega^+= X\setminus\overline{\Omega}$ is $X$.
\end{thm}

Note also that Theorem \ref{main}
easily implies the classical Hartogs-Bochner theorem
for bounded domains, since holomorphic extension from the complement of
a closed round ball to $\C^n$ (containing $\overline{\Omega}$)
can be proved by combining the one-dimensional Cauchy formula
along parallel slices.

Since pseudoconvex domains coincide with their envelope of holomorphy,
we immediately obtain

\begin{cor}\label{com-hyp}
If $\Omega$ and $M$ are as in Theorem \ref{main} and $\Omega^+$ is contained
in a pseudoconvex \emph{proper} subdomain of $\C^n$, then Hartogs extension fails.
This holds in particular if $\overline{\Omega}$ contains a closed complex
subvariety of $\C^n$ of dimension $n-1$.
\end{cor}

If $M$ is unbounded the assumptions of Theorem \ref{main}
are symmetric with respect to the sides $\Omega^\pm$.
In view of Corollary \ref{com-hyp}, the reader may wonder
about similarities with the theorem of Tr\'epreau \cite{trep1} 
on local extension of CR functions from real hypersurfaces.
We will elaborate on the relation between the two results 
in Section \ref{sec:remove-sing}.

The picture changes significantly, if we restrict to extension of \emph{smooth}
CR functions from the boundary of $\Omega$. 
Obviously extension still fails if $\Omega$ contains a complex hypersurface.
The case where $\overline\Omega$ contains a complex hypersurface but $\Omega$ does not
is more delicate:
Examples constructed in \cite{BDS3} by the first two authors and Z.~S{\l}odkowski
show that simultaneous extension of smooth CR functions to $\Omega$
may be valid or not, depending on a finer geometry of intersection
of the complex hypersurface and the boundary.
Combining Theorem \ref{main} with \cite[Section 5]{BDS3} we get

\begin{cor}\label{smooth-distr}
There are domains $D\subset\C^2$ with smooth connected boundary
such that Hartogs extension fails but every $\cC^\infty$-smooth
CR function on $\partial D$ has a holomorphic extension to $D$
(which is smooth up to $\partial D$).
\end{cor}

For more results and questions on particular domains (bounded or unbounded)
we refer to
\cite{BDS1, BDS2, BDS3, BDS4, DM, HST1, HST2, L, Sar}.
It may also be interesting to study
domains obtained by intersecting smoothly bounded domains in $\C\P^2$
with $\C^2$, for example with respect to extension from parts of the boundary,
see \cite{NP} for domains in $\C^2$ and further references.

%Theorem \ref{main} includes domains $\Omega$ with strictly pseudoconvex boundary
%as studied in \cite{HST2}.
%\marginpar{\textcolor[rgb]{1.00,0.00,0.00}{Is it a correct citation?}}

%If the minimality assumption is dropped, there are additional subtleties.
%So it may occur that the envelope of $\Omega^+$ is
%$\C^n\setminus A$ where $A$ is a CR invariant proper subset of $M$.
%Then for every CR distribution $u$ on $M$ its restriction to
%$M\setminus A$ has a holomorphic extension to $\Omega$, but it
%may happen that this extension does not attain $u$ as distributional
%boundary values along $A$ (see Example \ref{remove}).
%At the moment of writing, we do not know whether a CR function
%without Hartogs extension always exists in this case.

The paper is organised as follows: After some preliminaries collected in Section 2,
we prove the easier direction in Theorem \ref{main} in Section 3.
More precisely, we show how properties of the envelope of holomorphy of $\Omega^+$
imply Hartogs extension by using jump formulas and $\overline{\partial}$-methods.
The converse direction is treated
in the Sections~\ref{secloc} and \ref{sec:obstructions}. 
Section~\ref{secloc} contains topological preparations, which permit
in particular, to localise to envelopes of thin collars of the domain. Section 5
completes the proof of Theorem \ref{main}. The main ingredient is the use of
holomorphic functions with polynomial growth in order to construct
nonextendible CR distributions. 
In Section~\ref{sec:stein}, a generalisation of the main theorem to domains 
in Stein manifolds is given.  
The final section relates our result
to the topic of removable singularities. 
More precisely, we analyse obstructions to extension
confined to $M$ and exhibit analogies to Tr\'epreau's theorem.

{\bf Acknowledgments:} The third author would like to thank the
School of Mathematical \& Statistical Science at
the Arizona State University
in Tempe for its hospitality during a visit, which allowed to start the present research.

\section{Preliminaries}

{\bf Riemann domains.}
First we recall basic material on Riemann domains and envelopes of holomorphy,
referring to the monograph \cite{JP} for a thorough introduction.
For a domain $D\subset\C^n$, we denote by
$\pi_D:{\textsf E}(D)\rightarrow\C^n$ its envelope of holomorphy.
It is a Riemann domain (i.e.~$\pi_D$ is a local biholomorphism)
and there is a canonical embedding
$\iota_D:D\hookrightarrow{\textsf E}(D)$ satisfying
$\pi_D\circ\iota_D=\mbox{id}_D$, allowing us to identify
$D$ with $\iota_D(D)\subset{\textsf E}(D)$.
A classical theorem based on the solution of the Levi problem
tells that ${\textsf E}(D)$ is Stein.

Following Grauert and Remmert \cite{GR},
one may associate to every Riemann domain $\pi:X\rightarrow\C^n$
an abstract closure $\overline{\pi}:\overline{X}\rightarrow\C^n$
and an abstract boundary $bX=\overline{X}\setminus X$. Referring to
\cite[Section 1.5]{JP} for a careful treatment of the subtle construction,
we record that $\overline{X}$ is equipped with a natural topology
which restricts to the standard topology on $X$. Moreover $\overline{X}$
is the closure of $X$ in $\overline{X}$,
and $\overline{\pi}$ is the continuous extension
of $\pi$. For \emph{$\cC^1$-smoothly} bounded domains $D\subset\C^n$,
the abstract closure coincides with the usual one, but for
rough boundaries the abstract boundary $bD$ may be multi-sheeted
above the standard boundary $\partial D$.

\medskip

{\bf Distributions.}
Recall some basic facts on distributions on a real manifold $M$.
We consider a covering $\{\omega_j\}$ of $M$ by coordinate neighborhoods
$\omega_j$, i.e.~open sets equipped with diffeomorphisms
$\kappa_j:\omega_j\rightarrow\tilde{\omega}_j\subset\R^m$, $m=\dim M$.
Following \cite[Section 6.3]{H} a \emph{distribution on $M$} is given by
such a covering together with distributions
$u_j\in\cD'(\tilde{\omega}_j)$ satisfying
\begin{equation}\label{distrib}
u_i[\varphi]=u_j\left[\,|J_{\kappa_{ij}}|\,\varphi\circ \kappa_{ij} \right]
\end{equation}
for every $\varphi\in\cD(\kappa_i(\omega_i\cap\omega_j))$,
where $\kappa_{ij}=\kappa_i\circ\kappa_j^{-1}$
and $J_{\kappa_{ij}}$ is the Jacobian determinant
of the transition map.
This definition of distributions on $M$ is natural in so far
that every function $g\in\cC(M)$ identifies with
the distribution $g_j:\varphi\mapsto\int (g\circ\kappa^{-1}_j)\,\varphi\,dx$,
$\varphi\in\cD(\tilde{\omega}_j)$, because of the transformation formula.
Note on the other hand that distributions do not canonically correspond
to elements of the dual space of $\cD(M)$. In the proof of
Proposition \ref{suff}, we shall see how to use metrics
to this end.

%\begin{rem}\rm
%In the literature, the above formalism is often avoided by working with
%dual space of $\cD(M)$.
%\end{rem}

\medskip

{\bf CR distributions.}
For a $\cC^\infty$-smooth hypersurface $M$ in $\C^n$,
we say that a distribution $u=\{u_j\}$ is a \emph{CR distribution}
if the $u_j$ satisfy the tangential Cauchy-Riemann equations
in the weak sense.
To a function $f\in\cO(\Omega)$ satisfying (\ref{polgrowth})
we associate weak boundary values in the following way:
Locally we can represent $M$ as a graph
\begin{equation}\label{graph}
y_n=h(z_1,\ldots,z_{n-1},x_n)=h(z',x_n),
\end{equation}
with $h\in\cC^{\infty}(\tilde{\omega})$ with
$\tilde{\omega}^{\rm open}\subset\C^{n-1}_{z'}\times\R_{x_n}$,
so that $\Omega$ lies on the side $\{y_n>h\}$.
Then the distributions $f_\epsilon\in\cD'({\tilde{\omega}})$ defined by
\[
f_\epsilon[\phi]=\int f(z',x_n+i(h(z',x_n)+\epsilon))\,
\varphi(z',x_n)\,dx_1\,dy_1\ldots dx_{n-1}\,dy_{n-1}\,dx_n
\]
tend to a CR distribution $f^*\in\cD'({\tilde{\omega}})$
for $\epsilon\downarrow 0$, see \cite[Theorem 7.2.6]{BER}.
We can select a cover of $M$ by open sets $\omega_j$ graphed
over $\tilde{\omega}_j^{\rm open}\subset\C^{n-1}\times\R$
and obtain distributions $f^*_j\in\cD'({\tilde{\omega}_j})$
as above.
By the Baouendi-Treves approximation theorem we can locally
approximate the $f^*_j$ by (pullbacks of) restrictions of entire functions
and derive that the $f^*_j$ satisfy (\ref{distrib}).

\medskip

{\bf Boundary values of holomorphic functions.}
For detailed information on CR functions, we refer to \cite{B}, \cite{MP}.
We will need the following fact:
\emph{Let $D$ be a full neighborhood of $M$ in $\C^n$ and let $f\in\cO(D\setminus M)$
have polynomial growth towards $M$ from both sides. If near every $z\in M$ the two local
boundary values of $f$ from opposite sides coincide, then $f$ extends holomorphically
through $M$}.

We sketch a proof based on the hypoanalytic wave front set $WF_{ha}(u)$,
which is defined for CR distributions $u$ on $\cC^\infty$-smooth embedded CR manifolds,
and refer to \cite{trev} for a thorough introduction and the basic structure theorems
we use in the sequel.
By definition $WF_{ha}(u)$ is a $\R_{>0}$-invariant subset of the pointed characteristic bundle.
More precisely, the characteristic bundle is the real line bundle
\[
H^0 M=\bigcup_{p\in M}\{\xi\in T_p^* M:\xi|_{H_p M}\equiv 0\}\subset T^* M,
\]
and $WF_{ha}(u)$ is a subset of $H^0 M$ minus the zero section.
Now the two local weak boundary values $f^-=f^+=f^*$ coincide.
The existence of each of the local extensions $f^\pm$ rules out
one side of the zero section in $H^0 M$ from $WF_{ha}(f^*)$.
Hence $WF_{ha}(f^*)$ is locally empty, whence $f^*$ extends holomorphically
to an ambient neighborhood.

\section{Holomorphic extension}

The following proposition yields sufficiency in Theorem \ref{main}.

\begin{prop}\label{suff}
Let $\Omega\subset\C^n$, $n\geq 2$, be a domain with smooth connected
boun\-dary $M$.
If the envelope of holomorphy of $\C^n\setminus\overline{\Omega}$
is $\C^n$, Hartogs extension is valid for $\Omega$.
\end{prop}

We will indicate how the proof can be pieced together from known techniques.
Straightforward modifications yield versions for varying degrees of regularity,
for example for continuous CR functions defined on a $\cC^1$-smooth boundary.

{\bf Proof:}
On $M$ we select a smooth Riemannian metric $\mu$ and fix the orientation induced
on $M$ as the boundary of $\Omega$. We may restrict to graph representations as in
(\ref{graph}) such that
\[
dx_1\wedge dy_1\wedge\ldots\wedge dx_{n-1}\wedge dy_{n-1}\wedge dx_n
\]
is positive and write $\mu$'s volume form $\sigma$ in the local coordinates $\kappa_j$ as
\[
\sigma=\sigma_j\,dx_1\wedge dy_1\wedge\ldots\wedge dx_{n-1}\wedge dy_{n-1}\wedge dx_n.
\]
Since $\sigma_i=J_{\kappa_{ji}}\sigma_j$ by the transformation formula,
the coordinate-wise defined products $\sigma_j u_j$ glue to an element
of $u_\mu\in\cD'(M)$, if $u=\{u_j\}$ is a distribution.

To every CR distribution $u=\{u_j\}$ on $M$, we may canonically associate a current
$T_u$ on $\C^n$ of bidegree $(0,1)$ in the following way:
For a smooth compactly supported $(2n-1)$-form $\psi\in\cD_{(2n-1)}(M)$
the function $\psi/\sigma$ (defined as the unique function $\tilde{\psi}$
satisfying $\psi=\tilde{\psi}\sigma$) has compact support. Hence
\[
T_{u,M}[\psi]=(u_\mu)[\psi/\sigma]
\]
is a $(2n-1)$-dimensional current $T_{u,M}$ on $M$.
Writing
$$\omega_{x,y'}=dx_1\wedge dy_1\wedge\ldots\wedge dx_{n-1}\wedge dy_{n-1}\wedge dx_n,$$
the equality
\[
(u_j)_{\mu_j} \left[\frac{\psi}{\sigma_j\omega_{x,y'}}\right]
=u_j [\psi/\omega_{x,y'}]
\]
holds locally. Thus $T_{u,M}$ is independent of $\mu$.
Finally we set
\[
T_u[\varphi]=T_{u,M}[(\iota_M)^*\varphi]
\]
for smooth $(n,n-1)$-forms $\varphi\in\cD_{(n,n-1)}(\C^n)$.
Here $\iota_M$ is the embedding of $M$ and $(\iota_M)^*\varphi$
is the pullback of $\varphi$ to $M$. Again
$T_u$ is $\mu$-independent and $\overline{\partial}$-closed.
For the last property, one may observe that $\overline{\partial}$-closedness
is a local property (since $M$ is properly embedded) and use the
Baouendi-Tr\`eves approximation theorem.

Since $H^1_{\overline{\partial}}(\C^n)=0$, the Dolbeault isomorphism
gives a distribution solution $f\in\cD'(\C^n)$ of
\[
\overline{\partial} f=T_u.
\]
Since $T_u$ has no mass outside $M$, $f$ restricts to holomorphic
functions $f^-$ on $\Omega^-=\Omega$ and $f^+$ on $\Omega^+=\C^n\setminus\overline{\Omega}$,
by elliptic regularity.
By \cite{C,LT}, $u$ is the jump from $f^-$ to $f^+$ in the following sense:
If $r\in\cC^{\infty}(U)$ is a local defining function of $M\cap U$, $U\Subset\C^n$,
then
\[
\int_{r=\epsilon}\! f^+ \varphi\; - \int_{r=-\epsilon}\! f^- \varphi \;
\longrightarrow\; T_u [\varphi],\;\;\mbox{if }\epsilon\downarrow 0,
\]
holds for all $\varphi\in\cC^{\infty}_{(n,n-1)}(U)$ with $\mbox{supp}\,\varphi\subset U$.
Now $f^+$ extends to a holomorphic function on $\C^n$ by assumption,
and
\[
\lim_{\epsilon\downarrow 0}\int_{r=\epsilon} f^+ \varphi
=\lim_{\epsilon\downarrow 0}\int_{r=-\epsilon} f^+ \varphi
=\int_M f^+ \varphi
\]
holds by continuity. Hence $f^+ - f^-$ defines the desired
extension of $u$ to $\Omega$. \qed

\section{Localization near $M$}\label{secloc}
\begin{figure}[h]
\parbox{.66\textwidth}{\resizebox{3.2in}{2.4in}{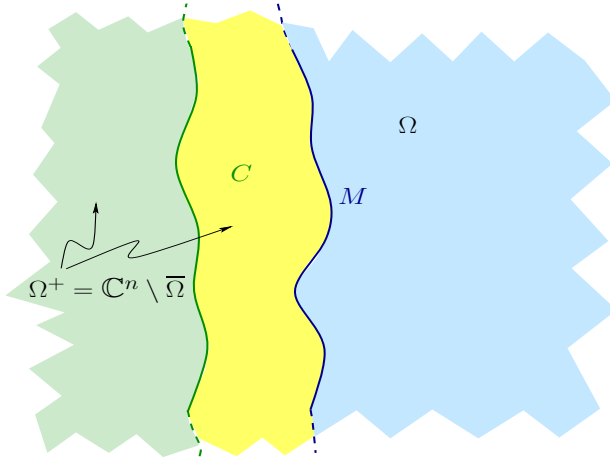}}
\parbox{.33\textwidth}
{\caption{\sl Outer collar $C$}
\label{fig:3-outer-collar}}
\end{figure}

A domain $C\subset\C^n\setminus\overline{\Omega}$
is called an \emph{outer collar of $M$} if $C\cup M$ is a
relative neighborhood of $M$ in $\C^n\setminus\Omega=\Omega^+\cup M$,
see Fig.~\ref{fig:3-outer-collar}.
Of course
$\Omega^+=\C^n\setminus\overline{\Omega}$ itself is an outer collar. In this section
we show that $\Omega^+$ can be replaced by an arbitrary outer collar in the assumptions
of Theorem \ref{main}. The most general version of our main result is

\begin{thm}\label{local}
For a domain $\Omega\subset\C^n$, $n\geq 2$, with connected smooth
boundary $M$ the following properties are equivalent:
\begin{itemize}
\item[\bf a)]
The envelope of holomorphy of $\Omega^+$ is $\C^n$.
\item[\bf b)]
For every outer collar $C$ of $M$, the canonical embedding
$\iota_{C}:C\hookrightarrow\textsf{E}(C)$ extends to a
(unique) lifting of $\Omega\cup M\cup C$ to $\textsf{E}(C)$.
\item[\bf c)]
There is an outer collar $C$ of $M$ such that $\iota_{C}$ extends
as in {\bf (b)}.
\item[\bf d)]
Every $u\in\cD'_{CR}(M)$ has a holomorphic extension to $\Omega$.
\end{itemize}
\end{thm}

Here we give the topological part of the proof, postponing extension
``\,{\bf (d)}"   to the next section.

{\bf Proof that {\bf (a)} $\Leftrightarrow$ {\bf (b)} $\Leftrightarrow$ {\bf (c)}:}
Since the implications {\bf (b)} $\Rightarrow$ {\bf (a)} $\Rightarrow$ {\bf (c)} are tautological,
it suffices to show {\bf (c)} $\Rightarrow$ {\bf (b)}.
We let $C_1$ be a collar as granted by {\bf (c)} and have to show the lifting property
for an arbitrary collar $C_2$.

\begin{lem}
If the lifting property holds for some subcollar $C'_2\subset C_2$, it also holds for $C_2$.
\end{lem}

{\bf Proof:}
The lifting property is equivalent to the fact that all $f\in\cO(C'_2)$ extend to
$\Omega\cup M\cup C'_2$. Applying this extension property to restrictions $g|_{C'_2}$,
$g\in\cO(C_2)$, we get the extension property and thereby the lifting property
for $C_2$. \qed

Hence it suffices to prove the lifting property for an appropriate subcollar of $C_2$,
allowing us to assume that $C_2\subset C_1$.

%
%First we observe that the lifting property for $C_2$ is valid as soon it holds for some subcollar $C'_2$.
%Therefore we may assume that $C_2\subset C_1$.

\begin{figure}[h]
\resizebox{5.6in}{3.4in}{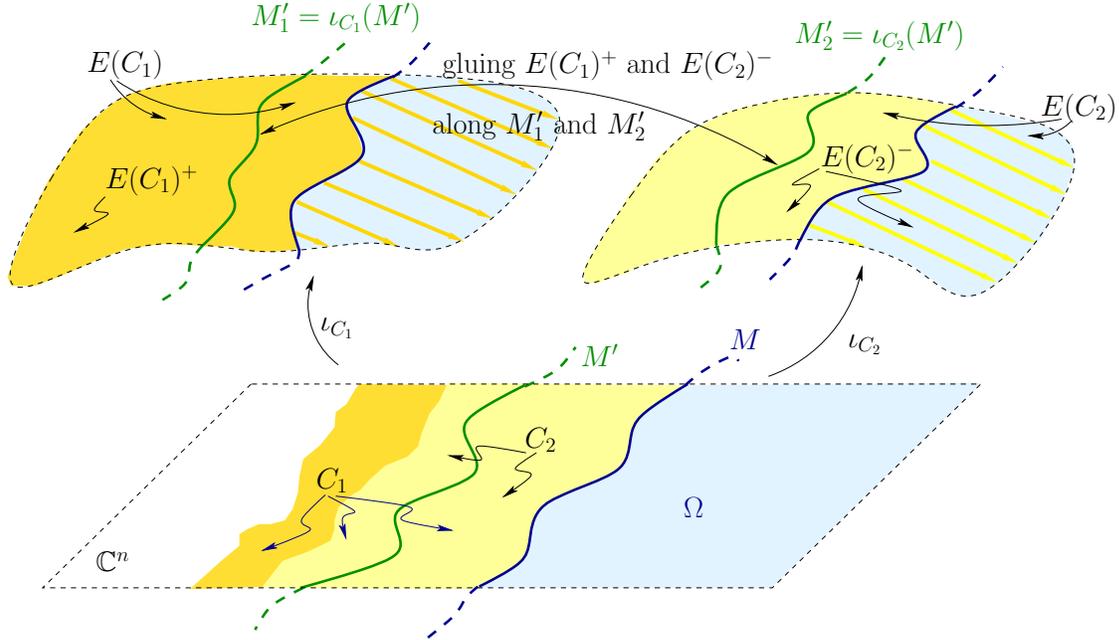}
\caption{\sl Two collars and corresponding envelopes}
\label{fig:4-two-collars}
\end{figure}

\begin{lem}\label{separate}
Let $M'$ be a smooth hypersurface obtained by isotoping $M$ into $C_2$.
Then $M'_1=\iota_{C_1}(M')$ disconnects ${\sf E}(C_1)$,
and $M'_2=\iota_{C_2}(M')$ disconnects ${\sf E}(C_2)$
\ (see Fig.~\ref{fig:4-two-collars}).
\end{lem}

{\bf Proof:} The argument is the same for $M'_1$ and $M'_2$.
Since $M'_1$ is connected, ${\sf E}(C_1)\setminus M'_1$ has at most two
connected components. If there is only one,
there is a smoothly embedded loop $\gamma\subset{\sf E}(C_1)$, which has exactly one
transverse intersection point with $M'_1$ (take a small arc transverse to $M'_1$
and link the endpoints by another arc that does not intersect $M'_1$).
By a result of Kerner \cite{K}, see also \cite{Stein}, $\iota_{C_1}$ induces
a surjective homomorphism $(\iota_{C_1})_*:\pi({C_1})\hookrightarrow\pi({\sf E}({C_1}))$
between the fundamental groups. Hence there is a loop $\tilde{\gamma}\subset {C_1}$
such that $\iota_{C_1}(\tilde{\gamma})$ and $\gamma$ are homotopic within ${\sf E}({C_1})$.

We use intersection numbers of oriented loops $\lambda\subset\C^n$
with $M'$, which can be defined as follows:
Let $\Omega'$ be the domain in $\C^n$ bounded by $M'$ and containing $\Omega$.
For $\lambda$ transverse to $M'$,
we compute the intersection number by
subtracting the number of points where $\lambda$ enters $\Omega'$
from the number of points where $\lambda$ leaves $\Omega'$.
The definition extends to general $\lambda$ because the intersection
number is homotopy invariant, see \cite{GP} for details.

Since $M\cap {C_1}=\emptyset$ the intersection number of $\tilde{\gamma}$ and $M$
is zero and the same holds for the intersection of $\tilde{\gamma}$ and $M^\prime$
(which is isotopic to $M$).
Pushing forwards by $\iota_{C_1}$, we get zero intersection
number between $\iota_{C_1}(\tilde{\gamma})$ and $M'_1$ (the intersection
number is calculated locally at the intersection points).
This contradicts the stability of intersection numbers under
homotopy and the fact that the intersection number of $\gamma$
and $M'_1$ is $\pm 1$. \qed

{\sl Continuation of the proof {\bf (a)} $\Leftrightarrow$ {\bf (b)} $\Leftrightarrow$ {\bf (c)}:}
%Let us now assume that the lifting property claimed in {\bf (b)} fails for $C_2$.
Denote by ${\sf E}(C_2)^-$ the connected component of ${\sf E}(C_2)\setminus M'_2$
which lies on the side of $\Omega$, see Fig.~\ref{fig:4-two-collars}.
More precisely, $M$ can be identified,
via the canonic embedding $C_2\hookrightarrow {\sf E}(C_2)$,
with a subset of the abstract closure of ${\sf E}(C_2)$, and ${\sf E}(C_2)^-$ is the
connected component containing $M$ in its closure. It suffices to show that
the domain $\Omega'$ considered in the proof of Lemma \ref{separate} lifts
to ${\sf E}(C_2)^-$ biholomorphically.

To this end, we construct a new Riemann domain $\pi_{X'}:X'\rightarrow\C^n$
by gluing ${\sf E}(C_2)^-$ with ${\sf E}(C_1)^+$ along $M'$.
More precisely, ${\sf E}(C_1)^+$ is the connected component of ${\sf E}(C_1)\setminus M'_1$
on the side opposite to $M$, and the gluing identifies $M'_2$ with $M'_1$.
Denote by $\widetilde{M}$ the corresponding hypersurface of $X'$. Note that
there is a natural embedding $\iota'$ of $C_1$ into $X'$, which coincides
with $\iota_{C_1}$ along $M'$.

Since pseudoconvexity is a local property at points of the
abstract boundary and $X'$ is obtained by gluing two pseudoconvex
Riemann domains along a set in the interior, $X'$ is pseudoconvex
by \cite{DG}.
Treating the sides of $\widetilde{M}$
separately, we see that every $f\in\cO(C_1)$ extends to $X'$.
Hence $X'$ is an extension of $C_1$ in the terminology of \cite[Section 1.4]{JP},
and by the pseudoconvexity of $X'$ this extension is maximal.
Hence $X'$ is equivalent to the envelope of holomorphy ${\sf E}(C_1)$ as a Riemann
domain over $\C^n$, meaning that $\Omega'$ lifts to $X'$.
Since the lifting has image in the side of $\widetilde{M}$
which is equivalent with ${\sf E}(C_2)^-$, we have proved
the existence of the desired lifting for $C_2$.

This finishes the proof that the first three properties are
equivalent. The link to {\bf (d)} will be completed in the
subsequent section.
\qed

%\begin{rem}\rm
%Observe that Theorem \ref{local} combine
%to a proof of the classical Hartogs-Bochner theorem.
%To this end, we use the known fact that smooth boundaries
%of bounded domains  $\Omega$ are globally minimal
%and observe that ${\textsf E}(\Omega^-)=\C^n$
%is trivial for bounded $\Omega$ (extend from the complement
%of a sufficiently large ball). \qed
%\end{rem}

\section{Obstructions to Hartogs extension}\label{sec:obstructions}

In this section, we will prove the harder direction in Theorem \ref{main}.

\paragraph{Geometry of ${\textsf E}(\Omega^+)$:}
Recall that we consider a domain $\Omega=\Omega^-$
with smooth connected boundary $M$, and
assume that the envelope of holomorphy
$\pi_{\Omega^+}:{\textsf E}(\Omega^+)\rightarrow\C^n$
of $\Omega^+=\C^n\setminus\overline{\Omega}$
differs from $\C^n$. For notational simplicity we write
$X^+$ instead of ${\textsf E}(\Omega^+)$, $\pi$
instead of $\pi_{\Omega^+}$,
$\overline{\pi}$ instead of the continuous extension
$\overline{(\pi_{\Omega^+})}:\overline{X^+}\rightarrow\C^n$,
and $\iota:\Omega^+\hookrightarrow X^+$ instead of
$\iota_{\Omega^+}$.

\begin{figure}[h]
\parbox{.56\textwidth}{\resizebox{3.4in}{3.4in}{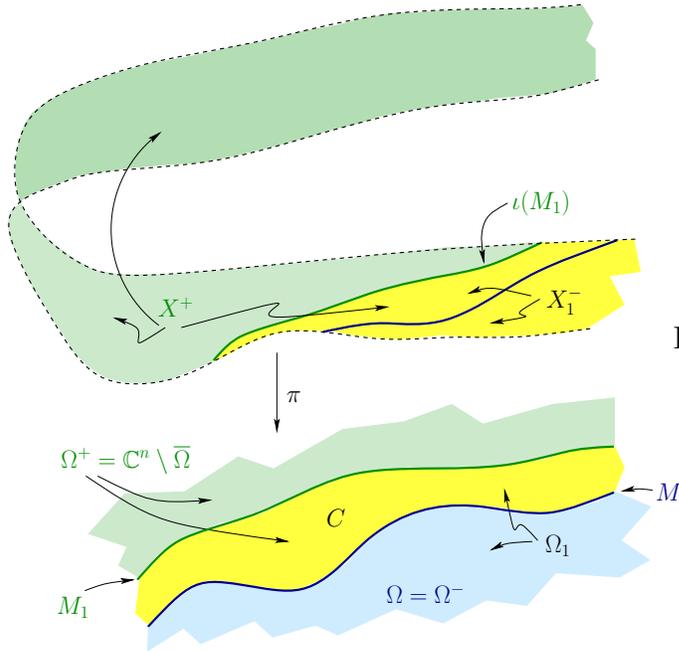}}
\parbox{.43\textwidth}
{\caption{\sl Multi-sheetedness}
\label{fig:5-multisheet}}
\end{figure}

Let $M_1$ be a smooth hypersurface obtained by slightly deforming
$M$ into $\Omega^+$. Let $\Omega_1$
be the domain that is bounded by $M_1$ and contains $\Omega$.
Then the intersection $C=\Omega_1\cap\Omega^+$ is a one-sided
collar of $M$, lying opposite to $\Omega$.
As in the proof of Theorem \ref{local}, we see that
$\iota(M_1)$ disconnects $X^+$ into two domains. Let $X^-_1$
be the connected component of $X^+\setminus\iota(M_1)$
that contains $C$.
Obviously $X^-_1=C\cup (X^+\setminus\iota(\Omega^+))$.
The arguments of Section \ref{secloc} actually 
allow us to identify $X^-_1$ with a subdomain of the envelope of $C$,
but we will not need this here. Note that $\pi(X^-_1)$
need not be contained in $\Omega_1$ and that there may
be multi-sheetedness over both sides of $M$,
as indicated in Figure \ref{fig:5-multisheet}.
%\begin{figure}[h]
%\parbox{.56\textwidth}{\resizebox{3.4in}{3.4in}{\input{5-multisheet.pstex_t}}}
%\parbox{.43\textwidth}
%{\caption{\sl Multi-sheetedness}
%\label{fig:5-multisheet}}
%\end{figure}

We distinguish two cases:

%Let $M_t$, $0\leq t\leq 1$, be a 1-parameter family of
%disjoint hypersurfaces, obtained by slightly deforming
%$M=M_0$ into $\Omega^-=\C^n\setminus\overline{\Omega}$.
%As in the proof of Proposition \ref{local} we see that
%each $M_t$, $t>0$, disconnects $X$ and bounds a connected
%component $X^-_t$ of $X\setminus X^-_t$
%on the side of $M$. Clearly the $X^-_t$
%decrease for $t\downarrow 0$ and
%$X^-=\bigcap_{0<t\leq 1} X^-_t$ is an open subset of $X$.
%Note that $X^-$ may be disconnected or empty and that
%$\pi(X^-)$ need not be contained in $\Omega$.

{\bf Case 1: $X^+$ is univalent.}
In this case we may identify both $X^+$ and  $X^-_1$
with domains in $\C^n$. Since $\Omega_1$ is one of the sides
of $M_1$ in $\C^n$ and $X^-_1$ contains $C$, $X^-_1$ is
contained in $\Omega_1$. It has to be a \emph{proper}
subset, for $X^+$ would be biholomorphic to $\C^n$ otherwise.
It follows that the abstract boundary $b X^-_1$ is the disjoint union of $M_1$ with
a \emph{nonvoid} set $S$ satisfying $\overline{\pi}_{\Omega^+}(S)\subset\Omega$.
Note that $\overline{\pi}_{\Omega^+}$ may become multisheeted on $S$, which therefore
cannot be identified with a subset of $\C^n$.

%We will distinguish two subcases:
%\begin{itemize}
%  \item {\bf Case 1.a:} $X^-_1\supset\Omega_1\setminus M$,
%  \item {\bf Case 1.b:} There is $z_0\in\Omega\setminus X^-_1$.
%\end{itemize}
%Case 1.b means that $bX$ contains some point lying over $\Omega$.
%Case 1.a cannot occur if $M$ is globally minimal, a consequence of the following
%\begin{lem}
%In the Case 1.a, $A=\C^n\setminus X$ is a CR invariant proper subset of $M$.
%\end{lem}
%
%{\bf Proof:}
%Clearly $A$ cannot coincide with $M$ since $X$ is connected.
%If $A$ fails to be CR invariant, it contains a point $z_0$
%where the local CR orbit is an open germ.
%By Tr\'epreau's theorem \cite{Trep}, there is
%one side of $M$ ($\Omega^+$ or $\Omega^-$)
%from which holomorphic functions extend to a neighborhood of
%$z_0$, in contradiction to pseudoconvexity of $X=\C^n\setminus A$. \qed

{\bf Case 2: $X^+$ is multi-sheeted}.
Obviously, this can only happen if $\iota(\Omega^+)$ is a proper subset of $X^+$.
For later use, we will only need that there is a $p_0\in X^+\setminus\iota(\Omega^+)$
such that the fiber $\pi^{-1}(\pi(p_0))$ contains at least two elements.

\paragraph{Construction of CR functions.}
Lifting the Euclidian distance, we get a Riemannian metric on any Riemann domain
$\pi:X\rightarrow\C^n$ and thereby the distance\footnote{
$\mbox{dist}(p,bX)$ is the supremum of all $r>0$ such that $B_r(\pi(p))$ can be
lifted to $X$ so that $\pi(p)$ is mapped to $p$}
$\mbox{dist}(p,bX)$ between a point $p\in X$
and the abstract boundary $bX$.
For a non\-ne\-gative integer $k$ we consider the Banach space
\[
\cO^{(k)}(X)=\{f\in\cO(X):f(p)\,\delta_X^{k}(p)\mbox{ is bounded on }X\},
\]
where 
\[ \delta_X(p) = \min\bigm(\mbox{dist}(p,bX),\frac{1}{\sqrt{1+|\pi(p)|^2}} \bigm) , \]
see \cite[\S 2.5]{JP} for detailed information. As observed in \cite{NP}, see also \cite{P2},
these spaces are useful for constructing CR distributions with
prescribed singularities.

\begin{lem}\label{f^*}
Let $f\in\cO^{(k)}(X^+)$ be given.

{\bf a)} The restriction $f^+=f|_{\Omega^+}$ has a unique CR distribution $f^*\in\cD'_{CR}(M)$
of order $k+1$ as weak boundary values on $M$.

{\bf b)} If $f^*$ has a holomorphic extension $f^-\in\cO(\Omega)$ then $f^-$ and $f^+$ glue
to an entire function. In particular $f^*$ is smooth.
\end{lem}

{\bf Proof:}
For $z\in\Omega^+$ close to $M$, we have
\[
|f_{\Omega^+}(z)|\leq C\,\mbox{dist}^{-k}(z,bX^+)\leq C\,\mbox{dist}^{-k}(z,M).
\]
Thus $f|_{\Omega^+}$ has at most polynomial growth towards $M$, and {\bf (a)} follows
from classical results on boundary values of holomorphic functions,
see \cite[Ch.~VII]{BER}.
If there is also an extension $f^-$ to the other side, the hypoanalytic wave front
set of $f^*$ is empty, meaning that $f^*$ is locally a restriction of a holomorphic
function, and {\bf (b)} follows from the uniqueness of holomorphic extension. \qed

From \cite{NP} we recall the following

\begin{lem}\label{sing}
Let $\pi_Y:Y\rightarrow\C^n$ be a pseudoconvex Riemann domain and $q$ a point
on the abstract boundary $bY$. Then there is a sequence
$Y\ni p_j\rightarrow q$ and a function $f\in\cO^{(2n+1)}(Y)$
such that $|f(p_j)|\rightarrow\infty$.
\end{lem}

Now we ready to construct a nonextendable CR distribution in the two
cases from our previous discussion of the geometry.
In {\bf Case 1}
there is a point $q\in bX^+\cap\pi_{X^+}^{-1}(\overline{\Omega})$,
a sequence $X^+\ni p_j\rightarrow q$ (convergence with respect to topology of the
abstract closure of $X^+$, which may be finer than the subspace topology
coming from $\C^n$) and a function $f\in\cO^{(2n+1)}(X^+)$ as in Lemma \ref{sing}.

We claim that the CR distribution $f^*$ associated to $f$ by Lemma \ref{f^*}
does not extend holomorphically to $\Omega$.
Otherwise Lemma \ref{f^*}, b) yields an entire function $F$.
The identity principle shows that $F$ is an extension
of $f$ (with $X^+$ considered as a subset of $\C^n$),
which is impossible since $|f(p_j)|\rightarrow\infty$.

In {\bf Case 2} we get a function $f\in\cO^{(6n+1)}(X)$ which separates the points
in the fiber $\pi^{-1}(\pi(p_0))$ from \cite[Proposition 2.5.5]{JP}.
Actually it suffices that $f$ attains different values at $p_0$ and a further
point $p_1\in\pi^{-1}(\pi(p_0))$.
Again we claim that the induced CR distribution $f^*$ does not extend to $\Omega$.
Otherwise the extension and $f|_{\Omega^+}$ glue along $M$ to an entire function $F$.
The identity principle yields $f=F\circ\pi$, and therefore $f(p_0)=f(p_1)$,
%Let $\gamma:[-1,1]\rightarrow X^+$ be a continuous path with $\gamma(-1)\in\Omega^+$,
%$\gamma(0)=p_0$ and $\gamma(1)=p_1$. Following simultaneously $\gamma$ and the
%projected path $\tilde{\gamma}=\pi\circ\gamma$, we see that
%$f(\gamma(0))=F(\tilde{\gamma}(0))=F(\tilde{\gamma}(1))=f(\gamma(1))$,
in contradiction to the choice of $f$.

\begin{rem}\rm
{\bf a)}
In some cases (for example if $M$ is strictly pseudoconvex at every point
and $\Omega$ is the pseudoconvex side)
the CR distributions we find to be obstructions to Hartogs extension
are smooth on $M$. In general this cannot always be achieved because
of examples constructed in \cite{BDS3}, see Corollary \ref{smooth-distr}.\\
{\bf b)}
It may happen that each CR distribution $u$ on $M$ possesses a weaker kind of
holomorphic extension to $\Omega$ which attains $u$ as weak boundary values
only along an open subset of $M$. We look at this in Section~\ref{sec:remove-sing}. \qed
\end{rem}

%For some sequence $p_j$ converging to the point $\tilde{z}_0\in\widetilde{M}$
%lying above $z_0$, Lemma \ref{sing} gives a function $f\in\cO^{(2n+1)}(X)$ with
%$|f(p_j)|\rightarrow\infty$. Its boundary value $f^*\in\cD'_{CR}(M)$ does not
%admit a holomorphic extension to $\Omega$. Otherwise $f$ would holomorphically
%extend to a neighborhood of $z_0$, in contradiction to the fact that it is
%unbounded near $z_0$.
%
%It continue with the Case 2.b$''$ where $\tilde{z}_0\in X$. There is another point
%$q_0\in\partial_X X^\ominus_1$ in the same fiber. By \cite{JP} there is
%a function $f\in\cO^{(6n+1)}(X)$ separating $\tilde{z}_0$ and $q_0$.
%Now $f^*$ is actually holomorphic near $z_0$. There cannot be a holomorphic
%extension to $\Omega$, as one can see by applying the identity principle
%along a smooth path which links $\tilde{z}_0$ and $q_0$ and remains in
%$X^\ominus_1$ (except at the extremal points).
%
%This settles the Case 2.b. Case 1.a is treated like 2.b$'$.
%Also Case 1.b can be treated similarly as Case 2.b'. Now we start with a
%function $f$ whose modulus explodes towards some $p_0\in bX$ above $\Omega$.
%Once again the identity principle applied along a suitable path shows that
%$f^*$ does not possess a holomorphic extension to $\Omega$.
%
%Case 2.a is handled in essentially the same way as 2.b$''$.
%Now we start from a function separating two points of
%$X^\ominus_1$ lying in the same fiber above some $z_0\in\Omega$.

{\bf Conclusion of the proof of Theorem \ref{main}:}
Sufficiency was shown in Section \ref{suff}. To derive necessity,
we argue by contraposition, assuming that the envelope of $\Omega^+$
differs from $\C^n$.
In each of the occuring cases, we have constructed a CR function
without extension to $\Omega$, which completes the proof of
Theorem \ref{main} and Theorem \ref{local}. \qed

\section{Generalisation of the main result to domains in Stein manifolds}
\label{sec:stein}

\begin{thm}\label{stein}
Theorem \ref{local} is still valid if $\C^n$ is replaced by a Stein manifold $Y$
of complex dimension $n\geq 2$.
\end{thm}

Most of the proof of Theorem \ref{local} is easily
generalised to Stein manifold.
The only ingredient 
which is specifically related to domains over $\C^n$ 
are the spaces $\cO^{k}(X)$. 
We shall give an extension to Stein manifolds, which is
less precise than the original results
but still sufficient for our needs.

Let $\pi:X\rightarrow Y$ be a Riemann domain over a complex manifold $Y$. 
Recall that an abstract boundary point $q\in bX$ 
can be specified by associating to every open neighborhood $U\subset Y$ 
of $p=\pi(q)$ the connected component $V$ of $\pi^{-1}(U)$ 
containing $q$ in its closure. 
If $U$ is relatively compact in $X$ 
and $z_1,\ldots,z_n$ are local holomorphic coordinates defined
in a neighborhood of $\overline{U}$, we may view $V$
as a Riemann domain over $\C^n$. 
Then we call $V$ \emph{adapted neighborhood}
and $z_1,\ldots,z_n$ \emph{adapted coordinates}.
We say that a function $f\in\cO(X)$ has 
\emph{polynomial growth of degree $k$ at $q\in bX$}, 
if there is an adapted neighborhood $V$ and $k\in\N_0$
such that $f|_V$ is an element of $\cO^{(k)}(V)$ with
respect to the corresponding adapted coordinates.
It is elementary to verify that this property does
not depend on the choice of adapted coordinates.
Define $\cO^{\rm pol}(X)$ as the algebra of 
\emph{holomorphic functions of polynomial growth},
i.e.~of all $f\in\cO(X)$ which have polynomial growth
at every $q\in bX$. 

In contrast to the Banach algebras $\cO^{(k)}(X)$, defined for $X$ spread over $\C^n$,
$\cO^{\rm pol}(Y)$ is only Fr\'echet in general. 
However, the following theorem is enough for constructing CR distributions
and permits to extend the proof of Theorem \ref{local} to Theorem \ref{stein}.

\begin{thm}\label{stein-pol}
Let $Y$ be a Stein manifold of dimension $n$ and $\pi:X\rightarrow Y$
a pseudoconvex Riemann domain over $Y$. Then the following hold:

{\bf a)} For every $q_0\in bX$, there is a sequence $x_j\in X$ with $x_j\rightarrow q_0$
and a function $f\in \cO^{\rm pol}(X)$ with $|f(x_j)|\rightarrow\infty$.

{\bf b)} For every pair $q_1,q_2\in X$ with $\pi(q_1)=\pi(q_2)$ and $q_1\not= q_2$,
there is  $f\in \cO^{\rm pol}(X)$ with $f(q_1)\not=f(q_2)$. 
\end{thm}

{\bf Proof:}
By the Bishop-Narasimhan-Remmert embedding theorem
(see \cite[Chap. VII, Sec.~C and Notes on p.~233]{GuRo}), we may assume
that $Y$ is a properly embedded complex submanifold of $\C^m$ for some $m\geq n$.
The normal bundle $\pi_N:N\rightarrow Y$ of $Y$ in $\C^m$
(i.e.~the bundle with fibers $T_y\C^m/T_y Y$, $y\in Y$)
has a natural holomorphic structure, with respect to which it is a Stein manifold.
We will identify $Y$ with the zero section of $N$. 
By \cite[Satz 3]{DG}, there is a holomorphic mapping 
$\Phi:N\rightarrow\C^m$, which is a locally biholomorphic 
outside some complex subvariety $A$ of $N\backslash Y$
of pure dimension $m-1$. 
Thus $\widetilde{\Phi}=\Phi|_{N\backslash A}:N\backslash A\rightarrow\C^m$ 
is a pseudoconvex Riemann domain.

Consider the Riemann domain $\tilde{\pi}:\widetilde{X}\rightarrow N\backslash A$
that is obtained gluing at every point $x\in X$ the corresponding fiber 
$\widetilde{\Phi}^{-1}(\pi (x))\backslash A$, or formally
$$\widetilde{X}=\{(x,v)\in X\times (N\backslash A):\widetilde{\Phi}(v)=\pi(x)\},\; \tilde{\pi}(x,v)=v.$$
It is straightforward to see that $\widetilde{X}$ is Stein 
and that $\alpha:x\mapsto (x,\pi(x))$ embeds $X$ into $\widetilde{X}$
in such a way that the image consists of the points 
lying above the zero section of $N$.
Moreover, $\Phi\circ\tilde{\pi}$ turns $\widetilde{X}$ into a Riemann domain over $\C^m$,
and we get the commutative diagram
\[
\begin{tikzcd}[row sep=large,column sep=huge]
       \widetilde X \arrow[r,"\tilde \pi"] & N\backslash A \arrow[d, "\widetilde \Phi"] \\
       X \arrow[r, hook] \arrow[u, hook, "\alpha"] & \mathbb C^m
\end{tikzcd}
\]
where the lower horizontal arrow denotes inclusion.

Hence we may apply results from \cite{JP} in the same way as in the proof of \cite[Lemma 2.2]{NP},
in order to obtain functions $f\in\cO^{(6m+1)}(\widetilde{X})$ which
explode at a given $\tilde{q}_0\in b\widetilde{X}$ or
separate two points $\tilde{q}_1\not=\tilde{q}_2$ 
lying in the same fiber of $\Phi\circ\tilde{\pi}$.
Here $\cO^{(6m+1)}(\widetilde{X})$ is defined with respect 
to the standard structure of $\C^m$. 

To prove (a) and (b), we choose $\tilde{q}_j=\alpha(q_j)$.
%the points lying over
%the zero section of $N$ and corresponding to $q,q_1,q_2$.
%Here $\cO^{6m+1}(\widetilde{X})$ is defined with respect 
%to the standard structure of $\C^m$. 
Since $\widetilde{X}$ is a product near every point of ${\tilde{\pi}^{-1}(Y)}$,
the restriction $f$ to $\alpha(X)$ 
defines an element of $\cO^{\rm pol}(X)$.
The proofs of Theorem \ref{stein-pol} and  
of Theorem \ref{stein} are complete. \qed

\section{CR orbits and removable singularities}\label{sec:remove-sing}

In this section we will slightly change our viewpoint by treating the 
sides of $M$ on equal footing.
Let $M\subset\C^n$ be a smooth real hypersurface
%, which we study locally at 
and $z\in M$. 
Consider ambient neighborhoods $\omega$ of $z$ 
such that $\omega\backslash H$ has exactly two connected components $\omega^\pm$. 
Then $z$ is called  
\emph{local obstruction point of $M$} if there is a neighborhood basis of $z$ 
by neighborhoods $\omega$ as above and functions $f^\pm\in\cO(\omega^\pm)$ 
which do not extend holomorphically to neighborhoods of $z$. 
Tr\'epreau's theorem says that $z$ is a local obstruction point of $M$
if and only if there is a local holomorphic hypersurface $Z$ satisfying 
$z\in Z\subset M$.
In the literature, it is customary to call $M$ \emph{minimal at $z$} 
iff $z$ is no local obstruction point.

%Here we elaborate the following subcase of Case 1 from the previous section:
For  $M$ as in  Theorem \ref{main},
with sides $\Omega^-=\Omega$ and $\Omega^+=\C^n\setminus\overline{\Omega}$ ,
let $\pi_\pm:X^\pm\rightarrow\C^n$ denote the envelopes
of holomorphy of $\Omega^{\pm}$,
$bX^\pm$ their abstract boundaries, 
$\overline{\pi}_\pm$ the continuous extensions
of $\pi_\pm$ 
to the abstract closures $X^\pm\cup bX^\pm$, 
and $\iota_{\pm}: \Omega^\pm\hookrightarrow X^\pm$ 
the canonical embeddings.
A point $z\in M$ is called a 
\emph{global obstruction point of $M$}
if there are liftings $z^-\in bX^-$ and $z_+\in bX^+$ satisfying
\begin{equation}\label{obs}
z^\pm\in bX^\pm\cap\overline{\iota_\pm(\Omega^\pm)}\cap \overline{\pi}_\pm^{\,-1}(z).
\end{equation}
Observe that (\ref{obs}) is equivalent to the existence of functions $f^\pm\in\cO(\Omega^\pm)$ without
holomorphic extension across $z$.  We study the \emph{global obstruction set} 
$M_{\textsf{obs}}\subset M$ of all global obstruction points of $M$.

To investigate the geometry of $M_{\textsf{obs}}$, we recall the notion of CR orbits.
Two points of $M$ lie in the same \emph{CR orbit}
if they can be linked by a piecewise smooth
CR curve, i.e.~a curve whose velocity vectors are contained in the complex
tangent bundle $H M=\bigcup_{p\in M}(T_p M\cap J_p T_p M)$. 
Obviously the CR orbits form a disjoint decomposition of $M$.
A fundamental result of Sussmann \cite{S} tells that every orbit is
an injectively immersed smooth manifold with real dimension
at least equal to the rank $2n-2$ of $HM$.
From basics facts on ordinary differential equations, 
it follows that orbits are either open subsets of $M$
or injectively immersed complex manifolds of complex dimension $n-1$,
%(also referred to as lower-dimensional orbits), 
see \cite[Section 3.1]{MP},
and that the union $M_{\textsf{hol}}$ of the lower-dimensional orbits
is closed in $M$. 

The following theorem shows that the the global obstruction set
is completely determined by the CR geometry of $M$.

%A set is called \emph{CR invariant} if it is the union of CR orbits.

%Of course, the set $M_{\textsf{obs}}\subset M$ only singles out 
%the obstructions for extension to $\Omega^\pm$ which lie on $M$.
%However it has remarkable geometric properties.

\begin{thm}\label{obs-set}
In the situation of  Theorem \ref{main},
the sets $M_{\textsf{obs}}$ and $M_{\textsf{hol}}$ coincide.
In particular, $M_{\textsf{obs}}\subset M$ is either empty or unbounded and of positive 
$(2n-2)$-dimensional Hausdorff measure.
\end{thm}

{\bf Proof:}
First we claim that $M_{\textsf{hol}}$ is contained in $M_{\textsf{obs}}$. 
To see this, we fix $z\in M_{\textsf{hol}}$.
%we consider first the case where $M_{\textsf{ld}}\not=M$ 
%and $\C^n\backslash M_{\textsf{ld}}$ is a domain. 
Since $\C^n\backslash M_{\textsf{hol}}$ is Stein,
there are functions $f^\pm\in\cO^{(6n+1)}(\C^n\backslash M_{\textsf{hol}})$
and sequences $\{z_j^\pm\}\subset\Omega^\pm$ approaching
$z$ from the $\Omega^\pm$-side, respectively, such that
$|f(z_j)|\rightarrow\infty$.
In $\overline{X^\pm}$, the lifted sequences $\{\iota_\pm(z_j^\pm)\}$ 
converge to elements $z^\pm\in bX^\pm$ above $z$, and the claim follows.

Observe that $M_{\textsf{obs}}$ is closed as the intersection of
the two closed sets
\[
\pi_\pm\left( bX^\pm\cap\overline{(\iota_\pm(\Omega^\pm))} \right).
\]
Fix $z_0\in M_{\textsf{obs}}$. For any $0<\epsilon\ll 1$ the open set 
$B_\epsilon(z_0)\backslash M$ has two connected components
$\omega^\pm\subset\Omega^\pm$. 
The envelope $\textsf{E}(\omega^\pm)$ is also a Riemann domain
over $X^\pm$ and it is readily verified that 
\[
z_0\in \overline{\pi}_{\omega^+}\!\left(b\textsf{E}(\omega^+)\cap\overline{\iota_{\omega^+}({\omega^+})}\right)
\cap\, \overline{\pi}_{\omega^-}\!\left(b\textsf{E}(\omega^-)\cap\overline{\iota_{\omega^-}({\omega^-})}\right).
\] 
Hence Tr\'epreau's theorem implies that there is a local complex hypersurface $Z\subset M$
passing through $z_0$. Since $Z$ is the CR orbit through $z_0$ with respect to a sufficiently small
neighborhood of $z_0$, $Z$ is smooth and tangent to $HM$.

We claim that a neighborhood of $z_0$ in $Z$ is contained in  $M_{\textsf{obs}}$.
Otherwise there is a point $z_1\in Z$ such that all functions in $\cO(\Omega^*)$,
where $*$ is one of the signs $+$ or $-$,
extend to a uniform ambient neighborhood of $z_1$. Now we get a contradiction to (\ref{obs})
from the general theorem about propagation of extension to full neighborhoods along
complex submanifolds of $M$. 
Below we provide some details on 
how to apply propagation arguments to envelopes of holomorphy. 

Consider the CR orbit $\cO(z_0,M)$ of $z_0$ in $M$. The proof of Theorem \ref{obs-set} 
will be complete, as soon as we have shown that $\cO(z_0,M)$ is a lower-dimensional orbit 
and satisfies
\begin{equation}\label{orb-obs}
\cO(z_0,M)\subset M_{\textsf{obs}}.
\end{equation}
Let us first show \eqref{orb-obs} in case that $\cO(z_0,M)$ is lower-dimensional. 
Then $\cO(z_0,M)$ can be parametrized by an injective
holomorphic immersion $\alpha:Z\hookrightarrow\cO(z_0,M)$ of a connected
$(n-1)$-dimensional complex manifold $Z$. Note that the
the manifold topology of $\cO(z_0,M)$, i.e.~the pushforward of the topology
of $Z$ under $\alpha$, may be finer than the topology induced from ambient space.
However, the above arguments imply that $M_{\textsf{obs}}\cap\cO(z_0,M)$
is both open and closed in $\cO(z_0,M)$ with respect to the manifold topology.
This proves \eqref{orb-obs} for $\cO(z_0,M)$ lower-dimensional.

It remains to rule out the case that $\cO(z_0,M)$ is open in $M$. 
Then $M$ is minimal at some point $z_1\in\cO(z_0,M)$ 
(otherwise $\cO(z_0,M)$ were foliated by complex hypersurfaces), 
and Tr\'epreau's theorem implies
that CR functions locally extend to one side of $M$. Since this property propagates
along CR orbits, CR functions extend to one side at every point of $\cO(z_0,M)$,
in particular at $z_0$.
Below we will outline how the information on extension of CR functions yields that $z_0$ 
is contained in the envelope of at least one of the domains $\Omega^\pm$. 
This contradicts $z_0\in M_{\textsf{obs}}$, and completes the proof of \eqref{orb-obs}.

Let us sketch the link between extension of CR functions from $M$ and the envelopes of $\Omega^\pm$. 
We will use the method of analytic discs, see \cite{BER,B,MP} for detailed information.
The tools necessary to realise the following outline are explained in the Chapters 4
and 5 of \cite{MP}, see also \cite{MP1} for more on deformation of discs.
An analytic disc is a mapping $A:\overline{\D}\rightarrow\C^n$ which is holomorphic in $\D$ 
and has some smoothness up to the boundary $\T=\partial\D$ 
(for our needs $\cC^{2,\alpha}$ with $0<\alpha<1$ is enough).
One works with discs attached to $M$ (i.e.~$A(\T)\subset M$)
or with boundaries close to $M$.
In our case, one starts from a chain of discs $A_j$, $j=1,\ldots,m$, attached to $M$
and linking $z_1$ and $z_0$ in the sense that $A_1(-1)=z_1$, $A_j(1)=A_{j+1}(-1)$,
$j=1,\ldots,m-1$, and $A_m(1)=z_0$. These discs are small in the sense that they
are attached to subsets of $M$ which can be represented as graphs and that the
local solution theory of the Bishop equation can be used to deform discs.
Since $z_1$ is a minimal point, we can sweep out one local side of $M$ 
at $z_1$ by images of a $1$-parameter family of discs. Then one uses
this open set $U_0$ attached to $M$ at $z_1$ in order to deform $A_1$
and produce a nearby disc $\tilde{A}_1$, whose image,
viewed as a parametrised surface, is transverse to $M$
at $\tilde{A}_1(1)=A_1(1)$. Sliding the $\tilde{A}_1$ in the directions transverse to 
$\frac{\partial\tilde{A}_1}{\partial\theta}$ (where $\T=\{e^{i\theta},\theta\in\R\}$)
yields a family that sweeps out a one-sided neighborhood $U_1$ attached to $M$
at $A_1(1)$. Note that $U_0$ and $U_1$ may lie on opposite
sides of $M$.

Iterating this procedure, we finally obtain a one-sided neighborhood $U_m$
attached at $A_m(1)$. The continuity principle applied to the underlying families of discs
shows that holomorphic functions defined in an \emph{arbitrarily thin} ambient neighborhood
of (a sufficiently large subset of) $M$
extend to the open sets $U_j$. By construction $U_m$ intersects one of the sides of $M$.
To fix ideas, we assume that this side is $\Omega^+$. Let $M_t$, $t\in[0,1]$, be a smooth 
$1$-parameter deformation of $M=M_0$ such that $M_t\subset\Omega^-$, $0<t\leq 1$.
Together the deformations of discs constructed above depend on finitely many parameters,
the dependence being $\cC^{2,\beta}$-smooth for some $\beta\in(0,\alpha)$.
Inspection of the Bishop equation shows that we may extend these deformations
to the parameter $t$ for $0\leq t\ll 1$. More precisely, we locally write the $M_t$
as families of graphs and obtain the families for $t>0$ by using the same data
as for $t=0$. Since the resulting discs depend 
$\cC^{2,\beta}$-smoothly on all parameters including $t$,
we get slightly deformed open sets $U_{m,t}$ attached to $M_t$ such that 
functions holomorphic near $M_t$ extend to $U_{m,t}$. If $t$ is sufficiently
close to $0$, we get $z_0\in U_{m,t}$, and hence that holomorphic
functions extend from $\Omega^-$ to a uniform neighborhood of $z_0$,
the desired contradiction. The proof of Theorem \ref{obs-set} is complete. \qed

As an application, we revisit a special case  in Theorem \ref{main}.

\begin{prop}\label{remov}
Let $\Omega\subset\C^n$, $n\geq 2$, be a domain with connected smooth
boundary $M$. Assume that $X^+={\textsf E}(\Omega^+)$ is univalent
and that $\C^n\setminus\overline{X^+}$ is contained in $M$.
Then we have $\C^n\setminus\overline{X^+}=M_\textsf{hol}$.
Moreover, $M_\textsf{hol}$ is \emph{removable} in the following sense:
For every CR distribution $u\in\cD'_{CR}(M\setminus A)$,
there is a function $\tilde{u}\in\cO(\Omega)$ which attains
$u$ as weak boundary value along $M\setminus A$.
\end{prop}

{\bf Proof:}
Clearly $M_\textsf{hol}$ is a proper subset of $M$, since otherwise
$\Omega^+$ would be Stein and coincide with $X^+$.
Theorem \ref{obs-set} directly implies that $M_\textsf{ld}$ and $X^+$ are disjoint.
If $z\in M\backslash M_\textsf{hol}$ we see like in the proof of Theorem \ref{obs-set}
that holomorphic functions extend through $z$ at least from on of the sides 
$\Omega^\pm$. Since $X^+$ is Stein and by assumption contains
both sides, we conclude $z\in X^+$, 
and the first part of the proposition follows.

As for removability,
we argue similarly as in the proof of Proposition \ref{suff}:
The complement of $M\backslash M_\textsf{hol}$ in $X^+=\C^n\backslash M_\textsf{hol}$
has two connected components $\Omega^\pm$.
Solving a suitable $\overline{\partial}$ equation
on the pseudoconvex domain $X^+$, we find $f^\pm\in\cO(\Omega^\pm)$
so that $u$ is the jump between $f^-$ and $f^+$ along $M\setminus A$,
in the sense of weak boundary values.
By assumption $f^+$ admits an extension $\tilde{f}^+\in\cO(X^+)$,
and $\tilde{u}=\tilde{f}^+|_\Omega-f^-$ is the desired extension
of $u$. \qed

We do not get more even if $u$ is a CR distribution on $M$
that is the global weak boundary value of a holomorphic function
on $\Omega^+$, as shown by

\begin{expl}\rm
Let $G=\{\zeta\in\C:\rho(\zeta)<0\}\Subset\C^1$ be a smoothly bounded disc such that 
$0\in\partial G$ and $T_0\,\partial G=\{\eta=\mbox{Re}(\zeta)=0\}$. In addition, we assume
that $G$ is strictly concave at $0$ and that $\nabla\rho(0)$ is proportional to $-\frac{\partial}{\partial\eta}$,
meaning that all $\zeta\in\partial G\backslash\{0\}$ close to $0$ are contained in $\{\eta<0\}$.
The unbounded domain
\[
\Omega=\left\{(z_1,z_2)\in\C^2:\rho\left(z_1 \exp \left(1+|z_2|^2\right)\right)<0 \right\}
\]
has smooth connected boundary $M$ homeomorphic to the cylinder $S^1\times\C$.
%The domain
%\[
%\Omega=\left\{(z_1,z_2)\in\C^2:\left|z_1-\exp \left(-(1+|z_2|^2)\right)\right|<\exp \left(-(1+|z_2|^2)\right)\right\}
%\]
%has connected boundary $M$ homeomorphic to the cylinder $S^1\times\C$.
It is routine to verify that $M$ decomposes into two CR orbits,
the $z_2$-axis $A$ and the open orbit $M\setminus A$.
Applying the continuity principle to families of complex lines parallel
to $A$ shows that $\C^2\setminus A$ is the envelope of
$\Omega^+=\C^2\setminus\overline{\Omega}$. Thus Proposition
\ref{remov} implies that every CR distribution defined on
$M\setminus A$ has a holomorphic extension to $\Omega$.
The function
\[
g(z)=\exp(1/z_1)|_{\Omega^+}
\]
is holomorphic and locally bounded along $\partial\Omega^+$.
In fact, $g$ is continuous near $\partial\Omega^+\backslash A$
and $|g|<1$ holds near every $z\in A$.
%(note that there are no uniform bounds along $A$,
%since $\Omega$ becomes thin for $|z_2|\gg 0$).
Hence its weak boundary value $g^*$ is a CR function
in $L^\infty_{loc}(M)$. Obviously the extension to $\Omega$ is
$g_\Omega=\exp(1/z_1)|_{\Omega}$. Note that $g_\Omega$ does not
have polynomial growth along $A$, meaning that $g^*$ is not the
weak boundary value of $g_\Omega$ along A. \qed
\end{expl}

\bigskip\bigskip

\noindent\small{School of Mathematical and Statistical Sciences,
Arizona State University, Tempe, AZ 85287, U.S.A.}

\vspace{-.05in}
\noindent {\it e-mail:} boggess@asu.edu

\medskip

\noindent  \small{Department of Mathematics,
Missouri University of Science and Technology,
Rolla, MO 65409, U.S.A.,  and
Faculty of Mathematics, Cardinal Stefan Wyszy\'nski University,
W\'oycickiego 1/3, \ 01-938 Warsaw, Poland}

\vspace{-.05in}
\noindent {\it e-mail:} romand@mst.edu

\medskip

\noindent\small{Department of Mathematics,
Mid-Sweden University, 85170 Sundsvall, Sweden, and
Institute of Mathematics, Jan Kochanowski University,
\'Swi\c etokrzyska 15,
\ 25-406 Kielce, Poland}

\vspace{-.05in}
\noindent {\it e-mail:} Egmont.Porten@miun.se

\end{document}